\documentclass[10pt]{amsart}
\usepackage{latexsym,amssymb,graphicx}

\newtheorem{thm}{Theorem}[section]
\newtheorem*{main}{Main Theorem}
\newtheorem{lem}[thm]{Lemma}
\newtheorem{cor}[thm]{Corollary}
\newtheorem{prop}[thm]{Proposition}
\theoremstyle{definition}
\newtheorem{defn}[thm]{Definition}
\newtheorem{quest}[thm]{Question}

\newtheorem{rem}[thm]{Remark}
\newtheorem{exmp}[thm]{Example}
\newtheorem*{ack}{Acknowledgments}

\newcommand{\Z}{\ensuremath{{\mathbb{Z}}}}
\newcommand{\R}{\ensuremath{{\mathbb{R}}}}

\newcommand{\D}{\ensuremath{{\mathcal{D}}}}
\newcommand{\A}{\ensuremath{{\mathcal{A}}}}
\newcommand{\W}{\ensuremath{{\mathcal{W}}}}
\newcommand{\cd}{\ensuremath{{\mathcal{C}}_{\D}}}
\newcommand{\khat}{\ensuremath{\widehat{K}_m}}
\newcommand{\xdhat}{\ensuremath{\widehat{\mathcal{X}}_{\D}}}
\newcommand{\xdelchat}{\ensuremath{\widehat{\mathcal{X}}_{\D_{\delta}}}}
\newcommand{\xd}{\ensuremath{\mathcal{X}_{\D}}}
\newcommand{\gpos}{G^{+}}

\newcommand{\lk}{{\mbox{Lk}}}
\newcommand{\ascend}{\lk_{\uparrow}}
\newcommand{\descend}{\lk_{\downarrow}}
\newcommand{\PD}{\ensuremath{\mathcal{PD}}}
\newcommand{\PDm}{\ensuremath{\PD_{\ngeq \mu}}}
\newcommand{\lf}{\ensuremath{\text{LF}}}
\newcommand{\rf}{\ensuremath{\text{RF}}}
\newcommand{\dnf}{\ensuremath{\text{DNF}}}

\begin{document}

\title[The homology of Garside groups]{Bestvina's normal
form  complex and the homology of Garside groups}

\author[R.~Charney]{Ruth Charney}
      \address{Department of Mathematics\\
               The Ohio State University\\
               231 W.~18th Ave\\
               Columbus, OH 43210}
      \email{charney@math.ohio-state.edu}

\author[J.~Meier]{John Meier}
      \address{Department of Mathematics\\
               Lafayette College\\
               Easton, PA 18042}
      \email{meierj@lafayette.edu}
      \thanks{Meier thanks The Ohio State University for hosting
      him while on leave from Lafayette. Charney thanks the
      Mathematical Institute at Oxford University for their hospitality.}

\author[K.~Whittlesey]{Kim Whittlesey}
      \address{Department of Mathematics\\
               University of Illinois at Urbana-Champaign\\
               1409 West Green Street\\
               Urbana, IL 61801}
      \email{kwhittle@math.uiuc.edu}

\keywords{}\date{\today}

\begin{abstract}
A Garside group is a group admitting a finite lattice generating
set $\D$.
Using techniques developed by  Bestvina for Artin groups of finite type,
we construct $K(\pi,1)$s for Garside groups.  This construction
shows that the (co)homology of any Garside group $G$ is
easily computed given the lattice $\D$, and there is a simple
sufficient condition that implies $G$ is a duality group.  The universal
covers of these $K(\pi,1)$s enjoy Bestvina's  weak non-positive curvature
condition.
Under a certain tameness condition, this implies that every
solvable subgroup of $G$ is virtually abelian.
\end{abstract}

\maketitle

\section{Introduction}\label{intro}

The main goal of this note is to establish the following result:

\begin{main}\label{thm:kgone}
Garside groups admit finite $K(\pi,1)$s.
\end{main}

A Garside group $G$ is the group of fractions
of a Garside monoid $\gpos$, where $\gpos$ contains
a `Garside element' $\Delta$ whose divisors form a finite
lattice $\D$ that generates $\gpos$.
(We give a complete definition in
\S\ref{sec:normalforms}.)   These groups, first introduced by
Dehornoy and Paris in \cite{dp}, generalize
Artin groups of finite type, and they share many formal
properties with Artin groups of finite type.

The construction of these  $K(\pi,1)$s is concrete and is based
on the lattice $\D$.  That is, given the lattice $\D$,
our proof of the Main Theorem is  constructive, and
one can easily compute
the (co)homology of the associated Garside group
(see Theorem~\ref{thm:resolution}).  Topological
properties of subcomplexes of the geometric realization
$|\D|$ give information on the end connectivity
of any Garside group, which yields a  simple
condition  implying that a given
Garside group is a duality group
(see \ref{cor:infinitygroups} and \ref{prop:duality}).

The construction of these $K(\pi,1)$s mimics
a construction given by Bestvina in
the case of Artin groups of finite type \cite{best}.
This construction is also used by Brady, as well as  Brady and Watt,
in their constructions of new $K(\pi,1)$s for Artin groups
of finite type  (\cite{brady} and \cite{bradywatt}).   The
argument given here indicates that this construction is ``functorial,''
in that it works for any group admitting a lattice generating
set, not just the Artin groups of finite type.  In particular
we extend Bestvina's results on a weak nonpositive
curvature condition for Artin groups of finite type to all
Garside groups.   Under a mild tameness condition on the Garside element
(Definition~\ref{defn:tame}),  this
implies that every solvable subgroup
of the Garside group is virtually a finitely generated
abelian group (Corollary~\ref{cor:solsubgrp}).

\begin{ack}
 We  thank Patrick Dehornoy
for conversations at geometric group theory conferences
in Montreal and Lyon, and in particular for pointing out
that he and Lafont have an independent method of computing
the homology of Garside groups \cite{dl}.
\end{ack}

\section{Garside groups and the Deligne normal forms}\label{sec:normalforms}

In terms of the word problem and normal forms, Garside
groups behave much like Artin groups of finite type \cite{smallg}.
Here we collect relevant information about the word problem
for Garside groups.  Most of this has to do with
issues in the positive monoid $\gpos$.

\medskip

Let $M$ be a monoid and let the \emph{indivisible} elements
(often called `atoms')
be those $m \in M$ such that $m \neq 1$ and if $m = ab$ then
either $a=1$ or $b=1$.  Let $||m||$ be the supremum  of the
lengths of all expressions of $m$ in terms of indivisible
elements.  The monoid $M$ is \emph{atomic} if it is
generated by its indivisible elements, and the norm $||m||$
of any element is finite.
In an atomic monoid one can define a \emph{partial order} via left
(or right) divisibility: $a < b$ if $a c = b$ for some $c \in M - \{1\}$;
$a \le b$ if $a c = b$ for some $c \in  M$.  Conversely, if $M$ is a
finitely generated
monoid where the partial order given by divisibility contains
no infinite descending chains, then $M$ is atomic, with $||a||$ being the
length of the longest expression.

\begin{defn}[Garside monoid]
An atomic monoid $M$ is \emph{Garside} if it satisfies
\begin{enumerate}
\item left and right cancellation laws hold in $M$,
\item any two elements of $M$ admit a least common multiple and
a greatest common divisor on both the left and the right,
\item there exists an element $\Delta$ such that the left
and right divisors of $\Delta$ are the same, there are
finitely many of them, and they form a set of generators
for $M$.
\end{enumerate}
The element $\Delta$ is called a \emph{Garside element}.
We denote the set of divisors of $\Delta$ by $\D$, and
call this generating set the \emph{simple divisors}.
\end{defn}

Garside monoids satisfy Ore's criterion, hence they
embed in their group of fractions, and thus
we may define a \emph{Garside group} to be the
group of fractions of a Garside monoid.

\begin{rem} Many of the properties of these groups were described in
the case of the braid groups by Garside in \cite{garside}; hence the
name ``Garside group''. When first introduced in \cite{dp}, however,
these groups were called ``small Gaussian'', while ``Garside'' was
used for a slightly more restrictive condition (which required that
$\Delta$ be the least common multiple of the indivisible elements).
 As it has developed that the definition
given above is more natural and more useful, the name ``Garside''
has now been generally adopted for this class of monoids and their
groups of fractions \cite{gargp}.
\end{rem}

The classic
examples of Garside groups are the Artin groups
of finite type, where the Garside element is commonly
denoted $\Delta$.  These groups also admit presentations
based on other monoids, such as the ones explored by Bessis,
Brady and Watt,
where the Garside element is not the usual $\Delta$, but
rather a Coxeter element $\delta$ (\cite{bessis}, \cite{brady},
and \cite{bradywatt}).
This collection  contains many other
groups (see \cite{smallg}, \cite{gargp}, \cite{garside}, \cite{dp} and
\cite{pic}).

We denote an Garside group by $G$ and let $\gpos$ be its
positive monoid. The inverse of an element $g \in G$ will be
denoted by $\overline{g}$.  We note that the partial order given
by left (or right) divisibility extends to the group
by defining $g < h$ when $\overline{g}h \in \gpos - \{1\}$, for
$g, h \in G$.  (Because Garside monoids are cancellative,
this partial order is equivalent to the partial order
defined above when restricted to the positive monoid.)

\medskip

Let $g \in \gpos$.  The \emph{left front} of $g$ is defined to
be the left gcd of $g$ and $\Delta$, $\lf(g) = g \wedge \Delta$.
If $\mu \in \D$  satisfies
$\mu < g$ then $\mu < \lf(g)$.  One can use left fronts to define
a normal form in $\gpos$, and this normal form is commonly
referred to as the \emph{left greedy normal form}.

\begin{prop}\label{prop:normalform}
\emph{(3.5 in \cite{smallg})}
Let $\gpos$ be  an Garside monoid.  Then $g$
may be uniquely represented as a product of simple
divisors, $g = \mu_{1} \mu_{2} \cdots \mu_{n}$, where
$\mu_{i} = \lf(\mu_{i} \cdots \mu_{n})$.
\end{prop}

\begin{prop}\label{prop:leftfront}
\emph{(3.10 in \cite{smallg})}
For all $g, h \in
\gpos$, $\lf(gh) = \lf(g \lf(h))$.
\end{prop}

There is an analogously defined \emph{right greedy} normal
form in which one begins by taking the (oddly named)
\emph{right front} $\rf(g)$ to be the right gcd of the
monoid element $g$ and the Garside element $\Delta$.
Then as in Proposition~\ref{prop:normalform} each $g \in \gpos$
can be uniquely represented as the product of simple divisors
$g = \mu_1 \cdots \mu_n$ where $\mu_i = \rf(\mu_1 \cdots \mu_i)$.
We make reference to this right greedy normal form in our
discussion of ascending and descending links in the
sections that follow.

\begin{prop}\label{prop:delta}
\emph{(2.2 and 2.3 in \cite{smallg})}
 Let $\gpos$ be an Garside monoid
with Garside element $\Delta$.  Then there is a permutation
$\sigma$ of the set of simple divisors $\D$ such that
\[
\Delta \mu = \sigma(\mu) \Delta
\]
for all $\mu \in \D$.  In particular, there is an $m$ such that
$\Delta^{m}$ is central.  Further, if $\mu \in \D$ can be expressed as
the product of $n$ indivisibles, then $\sigma(\mu)$ can also be
expressed as the product of $n$ indivisibles.
\end{prop}

\begin{defn}[Complements]
If $\mu \in \D$ then by definition there is an element $\mu^{\ast}
\in \gpos$ such that $\mu \mu^{\ast} = \Delta$.
  We call $\mu^{\ast}$ the \emph{right complement}
of $\mu$, and note that because
 $\mu^{\ast}$ is a right divisor of $\Delta$,
$\mu^{\ast} \in \D$, and thus $\D$ is closed under right
complements.

There is also a left complement of $\mu$, $^{\ast}\mu$, where $^{\ast}\mu
\mu = \Delta$.  If we right multiply this last equation by
$\mu^{\ast}$ we get $^{\ast}\mu \Delta = \Delta \mu^{\ast}$, hence
$^{\ast}\mu = \sigma(\mu^{\ast})$.  While this formula shows that
one does not have to use the notation $^\ast \mu$, we
do use this notation since it is easier to understand
than $\sigma(\mu^\ast)$.
\end{defn}

Since $\D$ is a generating set for the positive monoid $\gpos$,
and every element of $\D$ divides $\Delta$,
in order to represent elements of the group $G$, it suffices to invert the
Garside element.  Thus we get a set of normal forms for an Garside
group that closely parallel those of Deligne for Artin groups
of finite type.

\begin{thm}\label{thm:delignenf}
Every element $g \in G$ can be  expressed uniquely as
$g = \mu_{1} \cdots \mu_{k} \Delta^{n}$ ($n \in \Z$) where the prefix
$\mu_{1} \cdots \mu_{k}$ is in $\gpos$, it is in left greedy
normal form,
and $\mu_{1} \neq \Delta$.
\end{thm}

\begin{defn}[Deligne normal forms]
We refer to the normal form given in Theorem~\ref{thm:delignenf}
as the \emph{Deligne normal form} of the Garside group $G$.   We denote
the Deligne normal form of an element $g\in G$ by $\dnf(g)$. Thus
$\dnf$ can be thought of as a map from $G$ to the free monoid
$\{\D \cup \{\Delta^{-1}\}\}^{\ast}$, providing a (set theoretic)
section of the natural surjection $\{\D \cup \{\Delta^{-1}\}\}^{\ast}
\twoheadrightarrow G$.
\end{defn}

In the next section we'll be working with the Cayley graph of
a Garside group with respect to the set of simple divisors $\D$.  We'll
make use of the following lemma.

\begin{lem}\label{lem:onedelta}
Let $\mu_{1} \cdots \mu_{n}$ be a word in $\D^{\ast}$ in left greedy
normal form, where $\mu_{1} \neq \Delta$.  If $\eta \in \D$, then
the left greedy normal form for $\mu_{1} \cdots \mu_{n} \eta$ begins
with at most one $\Delta$.
\end{lem}

\begin{proof}
It is clear that right multiplication by $\eta$ can produce no $\Delta$
or one $\Delta$, so it suffices to show that $\Delta^{2} \nleq
\mu_{1} \cdots \mu_{n} \eta$.  If $\Delta^{2} \leq
\mu_{1} \cdots \mu_{n} \eta$, then $\Delta^{2} b = \mu_{1} \cdots \mu_{n}=
 \eta$
for some $b \in \gpos$.  Thus, shifting one of the $\Delta$s past
$b$ we get
\[
\Delta \sigma(b) \Delta = \mu_{1} \cdots \mu_{n} \eta,
\]
hence
\[
\Delta \sigma(b)\ ^{\ast}\eta = \mu_{1} \cdots \mu_{n}
\]
and therefore $\Delta \leq \mu_{1} \cdots \mu_{n}$, which
is a contradiction.
\end{proof}

\section{The $K(\pi,1)$ Complexes}\label{sec:complexes}

Recall that a \emph{flag complex} $X$ is a simplicial complex
where every complete subgraph on $n$-vertices in $X^{(1)}$
is the $1$-skeleton of an $(n-1)$-simplex
in $X$.  Thus a flag complex is determined by its $1$-skeleton.

Fix an Garside group $G$, whose positive monoid is $\gpos$, with
 Garside element
$\Delta$, and simple divisors  $\D$.
Let $\cd$ be the Cayley graph with respect to the
lattice generating set $\D$.  That is, $\cd$ has vertices
corresponding to the elements of $G$ and directed,
labelled edges corresponding to right multiplication
by elements of $\D$.  Let $\xdhat$ be the flag complex
induced by the graph $\cd$; thus each complete subgraph
of $\cd$ is the $1$-skeleton of a simplex in $\xdhat$.
The action of $G$ on $\cd$ by left multiplication extends
to an action of $G$ on $\xdhat$. This action is free since
it is free on the set of vertices and $G$ is torsion-free
\cite{torfree}. We establish our Main
Theorem by proving

\begin{thm}\label{thm:fullcomplex}
The complex $\xdhat$ is a contractible $G$-complex, and
$G\backslash \xdhat$ is a finite $K(G,1)$.
\end{thm}

In order to establish Theorem~\ref{thm:fullcomplex}
we show that $\xdhat$ admits a product structure.
Let $\xd$ be the full subcomplex of $\xdhat$ induced by vertices
associated with the set of all $g \in G$ where
the Deligne normal form of $g$ contains no $\Delta$s.
The complex $\xd$
can also  be thought of as a \emph{coset complex} where
the  vertices correspond to cosets of $\langle \Delta \rangle$. Coset
representatives can be taken from the positive monoid $\gpos$, and
represented by the prefixes that occur in the Deligne normal
forms.   There is an edge between two vertices in this complex
if their coset representatives differ by right multiplication
by some $\mu \in \D - \Delta$, and $\xd$ is the flag complex
induced by this graph.

\begin{lem}\label{lem:product}
The complex $\xdhat$ is homeomorphic to the
product $\xd \times \R$.
\end{lem}

\begin{proof}
A $k$-simplex in $\xd$ corresponds to  a collection of elements
$a_{i} \in \gpos$ where $\Delta \nless a_{i}$, $a_{0} < \cdots < a_{k}$
and $a_{k} < a_{0} \Delta$.  Similarly, a $k$-simplex in $\xdhat$
consists of elements $a_{i} \in G$ where $a_{0} < \cdots < a_{k}$
(using the extension of the partial order to $G$) and
$a_{k} \leq a_{0}\Delta$.  We give $\R$ the standard simplicial
structure where the vertices correspond to the integers.

There are projections from the $0$-skeleton of
$\xdhat$ onto  the $0$-skeleta of $\xd$ and $\R$ defined
via the Deligne normal forms.  Let $g \in G$ have
Deligne normal form $\dnf(g) = \mu_{0} \cdots \mu_{k} \Delta^{n}$.
Define $\pi_{+}$ to be the
map sending $g$ to $\mu_{0} \cdots \mu_{k}$ and $\pi_{\Delta}$
the map sending $g$ to $n$.  Both maps can be thought of as (set
theoretic) retractions of $\xdhat^{(0)}$; $\pi_{+}$ retracts the
vertices of $\xdhat$ onto the vertices of $\xd$ while $\pi_{\Delta}$
can be viewed as retracting the vertices of $\xdhat$ onto the
vertices corresponding to the infinite cyclic subgroup
$\langle \Delta \rangle < G$.

We form continuous maps from $\xdhat$ to $\xd$ and $\R$ by extending
$\pi_{+}$ and $\pi_{\Delta}$ linearly over simplices.
 Since $\xd$ is a flag complex, it
suffices to show that $\pi_{\Delta}$ and $\pi_{+}$ take edges
of $\xdhat$ to edges or vertices of their target spaces.
For the rest of this discussion, fix  an
edge $e$ in $\xdhat$ and let the bounding vertices of $e$
correspond to the
 group elements $a \Delta^{n}$ and $a \Delta^{n} \mu$ where $a \in \gpos$,
$\Delta \nless a$, and $\mu \in \D$.  The Deligne
normal form for $a \Delta^{n} \mu$
is then either $\dnf(a \sigma^{n}(\mu)) \Delta^{n}$, or $\dnf(b) \Delta^{n+1}$
if it's the case that $a = b(^{*}[\sigma^{n}(\mu)])$.
(The element $a \sigma^{n}(\mu)$
is divisible by at most one $\Delta$ by Lemma~\ref{lem:onedelta}.)

The map $\pi_{\Delta}$ takes $e$ of $\xdhat$ to the edge $[n, n+1]$
(when right multiplying by $\mu$ introduces a $\Delta$) or it collapses
$e$ to the  vertex corresponding to $n$ (when
$ ^{*}[\delta^{n}(\mu)]$ is not a right divisor of $a$).
Hence the map $\pi_{\Delta}$ extends to a map from $\xdhat$ onto
$\R$.

Similarly  the map $\pi_{+}$ extends to edges.  Let $e$ be
as before, so that $\pi_{+}$ maps the vertex $a \Delta^n$ to $a$
and
\[
\pi_{+}(a \Delta^{n} \mu) =
\begin{cases}
a \sigma^{n}(\mu) & \text{if $a$ is not right divisible by
$ ^{\ast}[\sigma^{n}(\mu)]$}\\
b = a \left( ^{\ast}[\sigma^{n}(\mu)]\right)^{-1} & \text{otherwise}
\end{cases}\ .
\]
So $\pi_{+}$ takes $e$ to an edge of $\xd$,
except in the case where  $\mu = \Delta$, in which case
$ ^{\ast}[\delta^{n}(\mu)] =  ^{\ast}\!\Delta = 1$, and $\pi_{+}$
collapses $e$ to a vertex.

Since $\pi_{+}$ and $\pi_{\Delta}$ extend to $\xdhat^{(1)}$,
and all three spaces are flag complexes, these maps extend to all
of $\xdhat$.  Thus we have a
continuous map $\pi: \xdhat \rightarrow \xd \times \R$ given
on the level of vertices by $\pi(g) = (\pi_{+}(g), \pi_{\Delta}(g))$.
It's clear that $\pi$ is surjective.

There is also a natural return map $\varpi: \xd \times \R
\rightarrow \xdhat$ that on the level of vertices is described by
$(a, n) \mapsto a \Delta^{n}$.  The cells in $\xd \times \R$
are of the form $\sigma \times I$, where
$\sigma$ is a $k$-simplex in $\xd$, and $I$ is a vertex or edge.
By definition the  simplex $\sigma$ corresponds to
a collection of elements
$a_{i} \in \gpos$ where $\Delta \nless a_{i}$, $a_{0} < \cdots < a_{k}$
and $a_{k} < a_{0} \Delta$.  If $I = [n,n+1]$ then
there is then a simplicial subdivision
of the cell $\sigma \times I$ induced by taking the simplices
\[
\left\{ \{a_{0} \Delta^{n} < \cdots < a_{k} \Delta^{n} < a_{0}
\Delta^{n+1}\}, \ldots , \{a_{k} \Delta^{n} < a_{0}\Delta^{n+1} <
 \cdots < a_{k} \Delta^{n+1}\}
\right\}\ .
\]
With this subdivision, $\xd \times \R$ is a simplicial complex
(in fact a flag complex),
and so in order to extend $\varpi$ to a map of complexes,
it suffices to check that
the  $\varpi$ can be extended to a simplicial map between
the $1$-skeleta.

In the simplicial subdivision of $\xd \times \R$,
there is an edge between $(a,n)$ and $(b,m)$ if $|m-n| \le 1$ and
$a \Delta^{n} < b \Delta^{m} \le a \Delta^{n+1}$,
 or similarly  $b \Delta^{m} < a \Delta^{n} \le b \Delta^{m+1}$.
Because the two cases are symmetric we may use the previous set of
inequalities
without loss of generality, and we note that either $m = n$ or $m = n+1$.

If $m=n$ then by dividing out $\Delta^{n}$ we see that
$a < b \le a \Delta$ hence $b= a \mu$ for some $\mu \in \D$.
Thus $(b,m)$ maps to the vertex associated with $a \mu \Delta^{n}
= a \Delta^{n} \sigma^{-n}(\mu)$ which is joined to
the vertex associated with $a \sigma^{n}$ by an edge labelled
$\sigma^{-n}(\mu)$.

If $m= n + 1$ then by dividing out $\Delta^{n}$ we see that
$a < b \Delta \le a \Delta$ hence $b = a \mu^{-1}$ for some
$\mu \in \D$.  Thus $b \Delta = a \mu^{\ast}$ and the
pair $(b,m)$ maps to the vertex associated with $a \mu^{\ast}
\Delta^{n} = a \Delta^{n}\sigma^{-n}(\mu^{\ast})$ and so there
is an edge joining $\varpi(a,n)$ and $\varpi(b,m)$ when
$(a,n)$ is joined to $(b,m)$ in the simplicial decomposition
of $\xd \times \R$.

Since the composition $\varpi \circ \pi$ is the
identity on $\xdhat$, the two complexes are homeomorphic.
\end{proof}

\begin{rem} The complex $\xd$, viewed as a coset complex,
 admits a left $G$ action. The kernel of this action is $\langle
\Delta^m \rangle$,
where $m$ is the smallest power of $\Delta$ which is central, and
hence it descends to an action by the quotient group
$G_{\Delta} = G/\langle \Delta^{m} \rangle$. The action of
$G_{\Delta}$ on $\xd$ has
finite stabilizers, namely conjugates of the subgroup
$\langle \Delta\rangle /\langle \Delta^m \rangle$.

The projection $\pi_+$ of $\xdhat$ onto $\xd$
is equivariant with respect to the $G$-actions, but the inclusion
of $\xd$ as a subcomplex of $\xdhat$ is not.
\end{rem}

\begin{lem}\label{lem:cosetcomplex}
The complex $\xd$ is contractible.
\end{lem}

\begin{proof}
We use the norm $||a||$ ---
where  $||a||$  is the maximum  length of a representative
of $a$ as a product of indivisible elements ---
to define  a Morse function $\nu: \xd \rightarrow [0,\infty)$.
On the level of vertices, which are represented by the
coset representatives, $\nu(a) = ||a||$;
the map is then  extended linearly
over the simplices.  To see that this map is
non-constant on edges, consider the two possibilities:
for $a \in \gpos$ with $\Delta \nleq a$, and $\mu \in \D - \{\Delta\}$,
the left greedy normal form of $a\mu$ begins with  zero or one
$\Delta$.  If it contains no $\Delta$, then $a\mu$ is another
coset representative and $\nu(a\mu) =  ||a\mu|| \ge ||a|| + ||\mu||
> ||a||$.  Otherwise $a = b\sigma(\mu^\ast)$ for some
coset representative $b \in \gpos$.
and right multiplication by $\mu$ corresponds to an edge from
$b$ to $a$.  Again, $\nu(a) =  ||b\sigma(\mu^\ast)||
\ge ||b|| + ||\sigma(\mu^\ast)|| > ||b||$.

In order to establish that $\xd$ is contractible it suffices
to show that the descending links of vertices are contractible
(cf. \cite{best}).  Identifying the vertices with elements
$a \in \gpos$ whose left greedy normal form does not begin
with $\Delta$, one can describe the descending link of $a$
as the subcomplex induced by the vertices
$\{b ~|~ ||b|| < ||a||
\mbox{ and }b \mu = a\mbox{ for some }\mu \in \D\}.$
Thus the vertices of  $\descend(a)$ correspond to those
$\mu \in \D$ where $\Delta \le a \mu$.   If we express
$a$ in \emph{right} greedy normal form, so
 $a = \mu_{1} \cdots \mu_{k}$ where $\mu_k = \rf(a)$ is the right
front of $a$, then $\Delta \le a \mu$ if and only if
$\Delta = \rf(a\mu)=\rf(\mu_k\mu)$. This occurs precisely when
$\mu_k^\ast \leq \mu$. Thus
the descending link is the subcomplex spanned by the
simple divisors $\mu \in \D -
\{\Delta\}$ that are greater than the right complement of the
right front of $a$, and the vertex corresponding to
$\rf(a)^{\ast}$ is a cone point for $\descend(a)$.
\end{proof}

Because $\xd$ is contractible, $\xdhat \simeq \xd \times \R$
is also contractible.  The $G$ action on $\xdhat$ is free,
and the quotient $G \backslash \xdhat$ is finite, hence
$G$ admits a finite $K(\pi,1)$.  In principle one can compute
the homology of $G$ from this $K(\pi,1)$.  This process can
be made quite concrete if one knows the lattice $\D$.

\begin{defn}
Let $\D_{n}$ denote the set of all ordered $n$-tuples of
elements of $\D$, $[\mu_{1} | \cdots |\mu_{n}]$, such that
the product $\mu_{1} \cdots \mu_{n}$ is an element of $\D$.
Note that the largest $n$ for which this set is non-empty
is the maximal length of $\Delta$ with respect to the
indivisible elements, that is $||\Delta||$.
\end{defn}

\begin{thm}\label{thm:resolution}   Let $G$ be a Garside group
with simple divisors $\D$, and let $B_{\D}$ be the
quotient of $\xdhat$ under the action of $G$. Then $B_{\D}$
has one $k$-cell for every element in $\D_k$.  In particular, if
$d=||\Delta||$, then
there is a free resolution of $\Z$ as a trivial $\Z G$-module
of the form
\[
0 \rightarrow \Z G^{|\D_d|} \rightarrow \cdots
\rightarrow \Z G^{|\D_{2}|} \rightarrow \Z G^{|D|} \rightarrow \Z G
\xrightarrow{\epsilon} \Z \rightarrow 0\ .
\]
\end{thm}

\begin{proof}
One can describe the complex $B_{\D} = G\backslash \xdhat$
using a variation on the standard bar construction.  Since
the $1$-skeleton of $\xdhat$ is the Cayley graph of $G$ with
respect to $\D$, the complex
$B_{\D}$ contains a single vertex, and a loop for each $\mu \in \D$.
Inductively one proceeds by considering each factorization
$\mu = \mu_{1} \mu_{2} \cdots \mu_{k}$
where $\mu$ and each $\mu_{i}$ are in $\D$.  Such a factorization gives
rise to a $k$-simplex $[\mu_{1}| \mu_{2} | \cdots |\mu_{k}]$ whose
codimension $1$ faces attach to the $(k-1)$ cells of $B_{\D}^{(k-1)}$
corresponding to $[\mu_{2}| \cdots |\mu_{k}]$,
$[\mu_{1}| \cdots |\mu_{i}\mu_{i+1}| \cdots |\mu_{k}]$, and
$[\mu_{1}| \cdots |\mu_{k-1}]$.  Thus the resolution of
$\Z$ as a trivial $\Z G$-module induced by the universal cover $\xdhat$
resonates with the standard bar resolution of a group, only
instead of using the entire multiplication table as in the bar resolution,
one restricts to the partial multiplication
in the lattice $\D$.
\end{proof}

\begin{quest}
Given an Garside group $G$, is it always possible to find an Garside
monoid $M$ where $G$ is the group of fractions of $M$, and
the Garside element $\Delta$ satisfies $||\Delta|| = \text{cd}(G)$?
One suspects not, but this is possible for Artin groups of finite
type, as we show in \S\ref{sec:examples}.
\end{quest}

\section{The end connectivity of Garside
groups}\label{sec:infinity}

Just as one can compute the (co)homology of an Garside group
using the poset $\D$, one can also get
concrete information about the end connectivity of $\xd$
by looking at the (co)homology of certain
subcomplexes of the geometric realization of the poset
of simple divisors $|\D|$.  Recall that
a contractible, locally finite complex $X$ is \emph{$n$-connected
at infinity} if given any compact subcomplex $C \subset X$,
there is a subcomplex $D \subset X$ such that every map
$\phi: S^i \rightarrow X - D$ extends to a map $\widehat{\phi}:
B^{i+1} \rightarrow X - C$, for $-1 \le i \le n$.  The property
of being $0$-connected at infinity is usually called ``one
ended'' and $1$-connected at infinity is usually referred to
as ``simply connected at infinity''.  Similarly $X$ is
\emph{$n$-acyclic at infinity} if given any compact subcomplex $C \subset X$,
there is a subcomplex $D \subset X$ such that
$i$-cycles supported  outside of $D$ bound $(i+1)$-chains supported
outside of $C$.  If $G$ admits a finite $K(\pi,1)$, then the end
connectivity of the universal covers of all finite $K(\pi,1)$s is the
same, and the end connectivity  is
therefore a property of the group $G$.

\medskip

In \S\ref{sec:complexes} we  establish that $\xd$
is contractible by examining
the descending links of vertices with respect to the
Morse function $\nu: \xd \rightarrow [0,\infty)$.  The
ascending links of vertices  can
be used to describe the connectivity at infinity
of $\xd$.  Using this idea, Bestvina was able to show that
if $A$ is an Artin group of finite type, with cohomological
dimension $n$, then $A$ is $(n-2)$-connected at infinity \cite{best}.
(As we discuss below, this implies that $A$ is a duality group.)  Arbitrary
Garside groups are not highly connected at infinity, and in particular,
they are not all duality groups.   The
main result in this section is that the lattice $\D$
determines the connectivity at infinity of $\xd$ and
hence it determines the connectivity at infinity of the
associated Garside group.

\medskip

Let $\PD$ be the subposet $\D - \{\Delta\}$,
and for any element $\mu \in \D$ let $\PDm$ be
 the subposet consisting of all elements which do not have $\mu$ as
a left divisor
\[
\PDm :=  \{\eta \in \D~|~\mu \nleq \eta \}\ .
\]
The ascending links of vertices can be described
in terms of the geometric realizations of these
subposets.  Namely, if $a \in \gpos$, then just as $\descend(a)$ is the
subposet of elements greater than or equal to
 the right complement of the right front $\rf(a)^\ast$;
the ascending link is then $\ascend(a) \simeq |\PD_{\ngeq \rf(a)^\ast}|$.

\begin{thm}\label{thm:infinityxd}  Let $G$ be an Garside group
with simple divisors $\D$, and let $\xd$ be the
coset complex.    If the reduced
homology groups
$\widetilde{H}_i(\PD)$
and $\widetilde{H}_i(\PDm)$, for $i \le n$, are trivial,
then $\xd$ is $n$-acyclic at infinity.
\end{thm}

\begin{proof}
We quickly sketch this argument since it is similar to
the one used by Bestvina for Artin groups of finite type \cite{best}
and by Bestvina and Feighn for $\mbox{Out}(F_n)$ \cite{bf}.

One starts by listing the coset representatives $\{1 = a_0, a_1, \dots,
a_m , \dots \}$ such that $||a_i|| \le ||a_j||$ if $i < j$.
Define $K_m$ to be the subcomplex of $\xd$ induced by the
vertices corresponding to $\{a_0, \dots, a_m\}$ and
$K_{-1}=\emptyset$.    The complement
$\xd - K_{m-1}$ is then the complement $\xd - K_{m}$ with the
ascending link $\ascend(a_m)$ coned off.  (The vertex corresponding
to $a_m$ forms the cone point.)

For $i>0$ the ascending link of the vertex associated to  $a_i$ is
homeomorphic to the complex
$|\PD_{\ngeq \rf(a_i)^\ast}|$;  the ascending link of $a_0 = 1$
is  $|\PD|$.  Let $C(\ascend(a))$ denote the cone on the
ascending link.  Applying Mayer-Vietoris to the union
\[
\xd - K_{m-1} = \{\xd - K_{m}\} \cup_{\ascend(a_m)} C(\ascend(a_m))
\]
yields
\[
\cdots \rightarrow
\widetilde{H}_i(\ascend(a_m)) \rightarrow
\widetilde{H}_i(\xd - K_{m}) \rightarrow
\widetilde{H}_i(\xd - K_{m-1}) \rightarrow \cdots
\]
where the term corresponding to
$C(\ascend(a_m))$ has been removed since cones are contractible.
Because the  complex $\xd$ is contractible, and each
$\widetilde{H}_i(\ascend(a_m))$ is assumed to be trivial (for $i \le n$),
we see by induction that $\widetilde{H}_i(\xd - K_m) = 0$
for all $i \le n$ and all $m$.
\end{proof}

\begin{cor}\label{cor:infinitygroups}
The Garside group $G$ is $(n+1)$-connected at infinity if
the reduced homology groups
$\widetilde{H}_i(\PD)$
and $\widetilde{H}_i(\PDm)$, for $i \le n$, are trivial.
\end{cor}

\begin{proof}
Since $\xdhat$ is a free, cocompact $G$-space, the end
connectivity of $G$ is the same as the end connectivity
of $\xdhat$.

The complex $\xdhat$ can be filtered by subcomplexes
of the form $\khat = K_m \times [-m,m]$ where $K_m$ is the compact subcomplex
of $\xd$ described in the proof of Theorem~\ref{thm:infinityxd}.
The complement
$\xdhat - \khat$ decomposes into the
union of two contractible pieces
$U_m = \{\xd \times [-m,\infty)\} - \khat$ and
$V_m = \{\xd \times (-\infty,m]\} - \khat$.  The
intersection
$U_m \cap V_m$ deformation retracts onto $\xd - K_m$, hence by
Mayer-Vietoris
$\widetilde{H}^{i}(\xdhat - \khat) \simeq
\widetilde{H}^{i-1}(\xd - K_m)$.  It follows by
Theorem~\ref{thm:infinityxd} that if the homology of $|\PD|$ and each
$|\PDm|$ is trivial for $i \le n$, then
$\widetilde{H}^{i}(\xdhat - \khat) = 0$ for
$i \le n+1$ and all $m$.

Since  $U_m$ and $V_m$ are contractible, van Kampen's Theorem
implies that  $\xdhat - \khat$ is simply
connected, assuming that the complexes $|\PD|$ and $|\PDm|$
are connected.   Thus by the Hurewicz Theorem, the complements
$\xdhat - \khat$ are actually  $(n+1)$-connected,
and therefore $G$ is $(n+1)$-connected at infinity.
\end{proof}

\begin{defn}[Duality groups]
A group $G$ of type FP (eg, a group with a finite $K(\pi,1)$-space)
is a \emph{duality group} if it is
$n$-dimensional and $(n-2)$-acyclic at infinity.  Equivalently,
it is an $n$-dimensional duality group if its cohomology
with group ring coefficients is torsion free and concentrated
in dimension $n$.   The term `duality' is appropriate since
in an $n$-dimensional duality group there is a natural isomorphism
between the group's homology and cohomology giving
$H_i(G,M) \simeq H^{n-i}(G, H^n(G;\Z G) \otimes M)$ for
all $i$ and all $G$-modules $M$.   The
module $H^n(G,\Z G)$ is called the dualizing module. (See
\cite{brownbk} for further background on duality groups.)
\end{defn}

\begin{prop}\label{prop:duality}
Let $G$ be a Garside group with simple divisors $\D$.  If the
cohomology of $|\PD|$ and each $|\PDm|$ is torsion free
and concentrated in dimension $n-2$, then $G$ is an
$n$-dimensional duality group.
\end{prop}

\begin{proof}
The cohomology group $H^i(G,\Z G)$ is isomorphic to
the direct limit of the induced system $H^{i-1}(\xdhat -
\khat)$.  By the proof of
Corollary~\ref{cor:infinitygroups} we see that
\[
\widetilde{H}^{i}(G, \Z G) =
\lim_{\rightarrow} \left\{
\widetilde{H}^{i - 2}(\xd - K_m) \right\}=
\widetilde{H}_c^{i-1}(\xd)
\]
Given that the cohomology of each  complex $|\PD|$
and $|\PDm|$ is
concentrated in dimension $n-2$, it follows from the Mayer-Vietoris
sequence in cohomology that there is a
short exact sequence
\[ 0 \rightarrow
\widetilde{H}^{i-2}(\ascend(a_n)) \rightarrow
\widetilde{H}^{i-2}(\xd - K_{n}) \rightarrow
\widetilde{H}^{i-2}(\xd - K_{n-1}) \rightarrow 0\ .
\]
Since the ascending links are finite,
$\widetilde{H}^{n-2}(\ascend(a_m))$ is free abelian,
and a  quick induction
shows that
\[
\widetilde{H}^{i-2}(\xd - K_n) \simeq \bigoplus_{i=0}^n
\widetilde{H}^{i-2}(|\PD_{\ngeq \rf(a_{i})^{\ast}}|)
\]
and so taking the direct limit we get
\[
\widetilde{H}^{i}(G, \Z G) = \widetilde{H}^{i-1}_{c}(\xd)
= \bigoplus_{i=0}^{\infty}
\widetilde{H}^{i - 2}(|\PD_{\ngeq \rf(a_{i})^{\ast}}|)
\]
hence $H^i(G,\Z G)$ is torsion free and concentrated in
dimension $n$ if the reduced cohomology of $|\PD|$ and each
$|\PDm|$ is torsion free and concentrated in dimension $(n-2)$.
\end{proof}

\section{Examples}\label{sec:examples}

The standard examples of Garside groups are the Artin groups
of finite type which are associated to the finite Coxeter groups.
Recall that a Coxeter system $(\W,S)$ consists of a finite set
$S=\{s_1,s_3,\dots ,s_n\}$ and a group $\W$ with presentation
\[
\W=\langle s_1,s_2,\dots ,s_n \mid s_i^2=1,
(s_is_j)^{m(i,j)}=1 \rangle
\]
where $m(i,j) \in \{2,3.\dots ,\infty\}$.
The associated Artin group $\A$ is the group with presentation
\[
\A=\langle s_1,s_2,\dots ,s_n \mid prod(i,j)=prod(j,i) \rangle
\]
where $prod(i,j)$ is the alternating product $s_is_js_i\dots$
of length $m(i,j)$. The pair  $(\A,S)$ is called an Artin system.
If $\W$ is finite, $\A$ is said to be of \emph{finite type}.
The map $\A \to \W$ sending each $s_i$ to the generator of the
same name is a quotient homomorphism.

For finite type $\A$, the monoid of positive words defined by
the presentation above is Garside. The indivisible
elements are the generators $S$,
and the lattice $\D$ is in one to
one correspondence with the Coxeter quotient $\W$, with
$\Delta$ corresponding to the maximal length element of $\W$.

Our construction of $\xdhat$, and especially $\xd$, is exactly
that of Bestvina for finite type Artin groups based on these
monoids \cite{best}.  However, there are other natural Garside
monoids whose groups of fractions are the Artin groups of finite
type.  The most well-known of these is the Birman-Ko-Lee braid
monoid \cite{birman} which was shown by Bessis, Digne, and Michel to
be a Garside monoid \cite{bdm}. This has been generalized by Bessis, Brady,
and Brady and Watt, to all Artin groups
of finite type (\cite{bessis}, \cite{brady}, and \cite{bradywatt}).
We briefly describe these monoids.

Let $(\W,S)$ be a finite Coxeter system and let $\A$ be the associated
Artin group.
Let $R$ be the set of all reflections in $\W$, that is, all conjugates
of $s_1, \dots s_n$. For $w \in \W$, let $|w|_R$ denote the minimal
length of an expression for $w$ as a product of reflections. If
$uv=w$ in $W$ and $|u|_R+|v|_R=|w|_R$, we write $u \leq_l w$ and $v
\leq_r w$.
These define partial orderings on $\W$. Let $\delta \in \W$ be the product
$\delta = s_{1} \cdots s_{n}$. Then $\delta$ is called
a Coxeter element of $\W$. (The choice of ordering of the
generators does not matter; permuting the generators gives a conjugate
Coxeter element.) The set
\[
\D_{\delta} = \{ w \in \W \mid w \leq_l \delta \}=
\{ w \in \W \mid w \leq_r \delta \}
\]
contains $R$ and is a lattice with respect to both partial
orderings.  Define $M_{\delta}$ to be the monoid with generators
$\D_{\delta}$ and relations $uv=w$ whenever $uv=w$ in $\W$, and  $u
\leq_l w \leq_l \delta$.  This monoid, which
Bessis calls the \emph{dual monoid}
is small Gaussian and its group of fractions is isomorphic to $\A$.

\begin{prop}\label{prop:dual}
Let $(\A,S)$ be an Artin system of finite type and let $M_{\delta}$
be its associated dual monoid.  Then
the complex $\xdelchat$ is a $K(\A,1)$ whose dimension $n=\vert S \vert$ =
is
the cohomological dimension of $\A$.
\end{prop}

It is easy to see that the cohomological dimension of $\A$ at least
$n$ since $\A$ contains a free abelian group of rank $n$.

\begin{rem} In the case of $2$- and $3$-generator finite type
Artin groups, the complex $\xdelchat$ is the same as the
$K(\A,1)$-complexes constructed by Brady in \cite{brady2}.
Brady shows that in these dimensions, this complex can be given
a piecewise Euclidean CAT(0) metric. Recent work of Brady and  McCammond,
however, shows that the analogous metric in higher dimensions will
fail to be CAT(0) for at least some of the Artin groups \cite{mcc}.
While the question of whether Artin groups are CAT(0) is still
open, in the next section we are able to show that Garside groups
satisfy a weak version of non-positive curvature.
\end{rem}

\begin{exmp}[Dihedral type]
Let $\A$ be an Artin group of finite dihedral type, that is,
the associated Coxeter group $\W$ is the dihedral group
of order $2n$, $D_{n}$.  In this case we make take the
standard generators of $\A$ to be $s$ and $t$ where the image
of $s$ and $t$ in $D_{n}$ correspond to reflections meeting
at an angle of $\pi/n$.  The Coxeter element $\delta = s\cdot t$
is then a rotation through an angle of $2 \pi/n$.

There are $n$ distinct reflections in $\W$.
Viewed as elements of $\A$, the set of reflections is
\[
R = \{s, t, \overline{t}st, st\overline{s}, \overline{ts}tst,
sts\overline{ts}, \ldots \}
\]
and the set of divisors is simply $\D_{\delta} = R \cup \{\delta\}$.
The complex $\xdelchat$ is thus the product of an $n$-valent
tree and $\R$. The orbit space $\A \backslash \xdelchat$
has a single vertex,
a loop based at that vertex for each element of $\D_{\delta}$
and a $2$-cell for each element of $R$:
$\{[s|t], [t| \overline{t}st], [ \overline{t}st| \overline{ts}tst],
\ldots\}$.  The link
of the vertex in $\A \backslash \xdelchat$ is topologically
the suspension of $n$ points, one for each element of $R$.
The suspension points correspond to the two ends,
$\{\delta^{+}, \delta^{-}\}$, of the loop $\delta$,
 and the arc corresponding
to a reflection $r$ is subdivided by two internal vertices
$r^{-}$ and $^{\ast}r^{+}$.
\end{exmp}

\begin{exmp}[Braid groups]
Let $\A$ be the braid group on $n+1$ strings. Then the reflections
$R$ consist of braids $r_{i,j},~i<j$, which interchange the $i$th
and $j$th strings and leave the other strings fixed (with a fixed
convention on over- or under-crossing) and the element $\delta$
is the braid which crosses the first string over to the last
position and shifts all the other strings one place to the left.
The complex $\A \backslash \xdelchat$ is $n$-dimensional
and has one $n$-cell for each expression of $\delta$ as a
product of $n$-reflections.
\end{exmp}

\section{Curvature}\label{sec:curvature}

A CAT(0) group is a group which acts properly and cocompactly on a
CAT(0) metric space. The existence of such an action has many
algebraic consequences for the group (see \cite{bh}.)
It is not known whether all Artin groups---or even those of finite
type---are CAT(0) groups. (See \cite{brady2}, \cite{bradymcc},
and \cite{bradycrisp} for some partial results on this question.)
Bestvina showed that the Artin groups of finite type, modulo their
centers, satisfy a combinatorial convexity property
which was good enough to get some interesting corollaries.
Many of these arguments work in the context of Garside groups as well.

\begin{defn}[Geodesics]
We define a \emph{geodesic} in $\xd$ to be any path whose sequence of
edge labels corresponds to a Deligne normal form with no $\Delta$s.
If $\dnf(v) = b_1 b_2 \ldots b_n \Delta^k$ with $b_n \not= \Delta$,
then $b_1, b_2, \ldots, b_n$ are the edge labels for the geodesic from
$1$ to $v$. The geodesic from $v$ to $w$ is the left translation by $v$ of
the geodesic from $1$ to $v^{-1}w$.
\end{defn}

\begin{defn}[Distance]
We define a \emph{nonsymmetric} distance $d(v,w)$ from a vertex $v$ to
a vertex $w$ in $\xd$ as follows: if $\dnf(v^{-1}w) = b_1 b_2 \ldots
b_n\Delta^k$, where $b_n \not= \Delta$, then $d(v,w) = ||b_1 b_2
\ldots b_n||$.
\end{defn}

\begin{lem}\label{lem:inverses}
The geodesic from $v$ to $w$ is the inverse of the geodesic from
$w$ to $v$.
\end{lem}

\begin{proof}
An edge path is a geodesic if and only if every subpath of length two
is a geodesic.  If $a$ and $b$ are simple divisors, then the path
$a,1,b$ is a geodesic if and only if ${}^*a b$ is in Deligne normal
form.  We must show that ${}^*b a$ is also in Deligne normal form. If
not, then ${}^*b a = {}^*b s a'$, where $s, {}^*b s,$ and $a'$ are
also simple divisors. Then
$$ ^*a b = \Delta (^*b a)^{-1} \Delta =
 \Delta (^*b s a')^{-1} \Delta = {}^*(a')(^*b s)^* =
{}^*a s ({}^*b s)^*$$
and so ${}^*a b$ was not in Deligne normal form.
\end{proof}

\begin{prop}\label{prop:triangle}
We can orient the edges of $\xd$ in terms of the Morse function.  Then
for any two vertices $v$ and $w$ in $\xd$, all edges of the geodesic
from $v$ to $w$ whose orientation points toward $v$ occur at the
beginning of the path, and these are followed by edges oriented
toward $w$.
\end{prop}

\begin{proof}
Assume to the contrary that $x,y,z$ are consecutive vertices on the
geodesic from $v$ to $w$ and that $||x|| < ||y|| > ||z||$.  Thus there
are divisors $b$ and $c$ in $\D$ such that $xb = y$ and $yc = z\Delta$.
It follows that the first letter in the left greedy normal form for
$x(bc) = z \Delta$ is $\Delta$.  But then by
Proposition~\ref{prop:leftfront}, $\lf(xb) = \Delta$.  But $xb = y$
and $\Delta \nleq y$.
\end{proof}

\begin{cor}\label{cor:nonposcurv}
If $p$ is a vertex on the geodesic from $v$ to $w$ and $v \ne p \ne
w$, and $u$ is any vertex, then there is a strict inequality:
$$ d(u,p) < \max\{d(u,v),d(u,w)\}.$$
\end{cor}

\begin{defn}[Centers]
Let $T$ be a finite set of vertices in $\xd$.
For a vertex $v$ in $\xd$, let $r(T,v)= max_{t \in T}{d(t,v)}$.
The \emph {circumscribed radius} $r(T)$ is the smallest integer $r$
such that $r=r(T,v)$ for some vertex $v$. Any such vertex is called
a \emph{center} of $T$.
\end{defn}

\begin{prop}\label{prop:centers}
The set of centers of $T$ span a simplex in $\xd$.
\end{prop}

\begin{proof}
Since $\xd$ is a flag complex, it suffices to show that any two
centers are connected by a single edge. Suppose $v_1$ and $v_2$ are
centers and that the geodesic between them passes through a third
vertex $w$. Then for every $t \in T$,
\[
d(t,w) < \max\{ d(t,v_1),d(t,v_2)\} \leq r(T)
\]
which contradicts the minimality of $r(T)$.
\end{proof}

Let $m$ be the smallest power of $\Delta$ which is central in $G$ and
let $G_{\Delta} = G/\langle \Delta^m \rangle$.  Recall that for
$\mu \in \D$, $\sigma(\mu) = \Delta \mu \Delta^{-1} \in \D$.

\begin{cor}\label{cor:fixedpoint}
Let $H$ be a finite subgroup of $G_{\Delta}$.  Then
$H$ fixes the barycenter of some simplex in $\xd$, and therefore
there are only finitely many conjugacy classes of finite subgroups
in $G_{\Delta}$.  Moreover, up to conjugacy, the finite subgroups
of $G_{\Delta}$ are of one of two types:
\begin{enumerate}
\item cyclic generated by $\mu \Delta^j$, where $\mu \in \D$ satisfies
$\mu \sigma^j(\mu) \sigma^{2j}(\mu) \cdots \sigma^{(t-1)j}(\mu) =\Delta$
for some $t \leq ||\Delta||$, or
\item the direct product of a cyclic group of type (1) and
$\langle\Delta^k\rangle$ where $\Delta^k$ commutes with $\mu$.
\end{enumerate}
\end{cor}

\begin{proof}
The first statement follows immediately from
Proposition~\ref{prop:centers}. Suppose $H$ stabilizes the simplex
$1 <\mu_1 <\mu_2 < \dots <\mu_{t-1}$ and acts
transitively on the vertices.  Following
Bestvina's argument (\cite{best}, Lemma 4.6), one can show that $H$
preserves the cyclic order on the vertices induced by the linear
ordering above. Let $h \in H$ be an element
which takes the vertex $1$ to the vertex $\mu = \mu_1$, so $h=\mu
\Delta^j$ for some $j$. Then for $i=1,2,\dots t-1$,
$h^i = \mu \sigma^j(\mu) \sigma^{2j}(\mu)
\cdots \sigma^{(i-1)j}(\mu)\Delta^{ij}$ and the length ordering
implies that
$\mu \sigma^j(\mu) \sigma^{2j}(\mu) \cdots \sigma^{(i-1)j}(\mu)= \mu_i$.
(If the product on the left produced a $\Delta$ at any stage, then
some $\mu_i$ would have shorter length than the previous one.)
It follows that $h^t$, which
fixes the entire simplex, must equal $\Delta^{tj+1}$ which gives the
equation in (1) above. Finally, note that if $h' \in H$ takes the
vertex $1$ to the vertex $\mu_k$, then $h'h^{-k}$ fixes the entire
simplex and hence must be a power of $\Delta$ which commutes with $\mu$.
\end{proof}

\begin{defn}[Translation Length] For $g \in G$, let $|g|_\D$
denote the word length of $g$ with respect to the generating set $\D$.
The \emph{translation length} of $g$, denoted $\tau(g)$, is
defined to be $$\tau(g) = \lim_{n\to\infty}\frac{|g^n|_\D}{n}$$
\end{defn}

In order to use the convexity property of $\xd$, we will need a
relation between the word length $|g|_\D$ and the distance $d(*,g(*))$
in $\xd$, where $G$ is understood to act on $\xd$ via the projection
onto $G_{\Delta}$.

\begin{defn}\label{defn:tame} The Garside element $\Delta$ is
\emph{tame} if there exists a constant
$c$ such that $||\Delta^n || \leq cn$ for all $n$.
\end{defn}

Note that by definition, the norm $||~||$ on $G^+$ satisfies a reverse
triangle inequality: $||ab|| \geq ||a|| + ||b||$. Thus, $||\Delta^n ||
\geq n||\Delta||$, but it is not  obvious that $\Delta$ is tame.
On the other hand, we know of no Garside groups for which this
condition fails.

\begin{lem}\label{lem:deltalength}
Suppose $\Delta$ is tame.
Then for all $g \in G$, $|g|_{\D} \geq \frac{1}{c}\,d(*,g(*))$.
\end{lem}

\begin{proof} The normal form $g=\mu_1\mu_2\cdots \mu_k\Delta^j$
is an expression of minimal length for $g$ as a word in $\D$
(\cite{dp}), thus $|g|_{\D}=k+|j|$. Let $a=\mu_1\mu_2\cdots \mu_k \in
G^{+}$. Then there exists $b \in G^{+}$ with
$ab=\Delta^k$. It follows that $d(*,g(*))=||a|| \leq ||a|| + ||b|| \leq
||\Delta^k|| \leq ck \leq c\,|g|_{\D}$.
\end{proof}

\begin{prop}\label{prop:translation length} If $G$ is a Garside
group with tame $\Delta$, then
the set of translation lengths of elements of $G$ is bounded
away from zero.
\end{prop}

\begin{proof}
Let $g \in G$. We first consider the case where the image of $g$ in
$G_{\Delta}$ has finite order. By
Corollary~\ref{cor:fixedpoint}, there exists an $N$ such that the
order of every torsion element in $G_{\Delta}$ divides $N$. Thus,
$g^N=\Delta^{km}$ for some $k \neq 0$ (since $G$ is
torsion-free). It follows that $\tau(g)=\frac{|km|}{N} \geq \frac{m}{N}$.

Now suppose the image of $g$ in $G_{\Delta}$ has infinite order.  Let
$T_n$ be the set of vertices in $\xd$, $T=\{*,g(*),g^2(*), \dots,
g^{n-1}(*)\}$. Let $r_n$ be the circumscribed radius of $T_n$.  We
claim that for any $n$, $r_n < r_{n+d}$ where $d = ||\Delta ||$.
For suppose $r_n = r_{n+d}$. If $z$ is a center of $T_{n+d}$, then
it is also a center for any subset of $T_{n+d}$ whose circumscribed
radius is $r_n$. In particular, $z$ is a center for each of the sets
$T_n, g(T_n), g^2(T_n), \dots ,g^d(T_n)$. It follows that $z, g(z),
g^2(z), \dots ,g^d(z)$ are all centers of $g^d(T_n)$, hence they span
a simplex. But the dimension of the complex $\xd$ is $d-1$ so we must
have $g^k(z)=z$ for some $k \leq d$. Since vertex stabilizers in
$G_{\Delta}$ are finite, this contradicts the assumption that the
image of $g$ has infinite order.

Thus we have $r_d < r_{2d} < r_{3d} ....$ and so $r_{nd} \geq n$ for
all $n$. In particular, for some $0<i\leq d$, we must have
$d(*,g^{nd+i}(*)) \geq n$,
and hence by Lemma~\ref{lem:deltalength}, $|g^{nd+i}|_{\D} \geq
\frac{n}{c}$. It now follows that $\tau(g) \geq \frac{1}{cd}.$
\end{proof}

\begin{cor}\label{cor:solsubgrp}  If $G$ is a Garside
group with tame $\Delta$, then
every solvable subgroup of $G$ is virtually a finitely
generated abelian group.
\end{cor}

\begin{proof} The proof is the same as for Bestvina's Corollaries 4.2 and
4.4.
\end{proof}

\end{document}